# MONTE CARLO LIKELIHOOD INFERENCE FOR MISSING DATA MODELS

BY YUN JU SUNG AND CHARLES J. GEYER

*University of Washington and University of Minnesota*

We describe a Monte Carlo method to approximate the maximum likelihood estimate (MLE), when there are missing data and the observed data likelihood is not available in closed form. This method uses simulated missing data that are independent and identically distributed and independent of the observed data. Our Monte Carlo approximation to the MLE is a consistent and asymptotically normal estimate of the minimizer $\theta^*$ of the Kullback–Leibler information, as both Monte Carlo and observed data sample sizes go to infinity simultaneously. Plug-in estimates of the asymptotic variance are provided for constructing confidence regions for $\theta^*$. We give Logit–Normal generalized linear mixed model examples, calculated using an R package.

**1. Introduction.** Missing data [20] either arise naturally—data that might have been observed are missing—or are intentionally chosen—a model includes random variables that are not observable (called *latent* variables or *random effects*). A normal mixture model or a generalized linear mixed model (GLMM) is an example of the latter. In either case, a model is specified for the *complete data* $(x, y)$, where $x$ is *missing* and $y$ is *observed*, by their joint density $f_\theta(x, y)$, also called the *complete data likelihood* (when considered as a function of $\theta$). The maximum likelihood estimator (MLE) maximizes the marginal density $f_\theta(y)$, also called the *observed data likelihood* (when considered as a function of $\theta$). This marginal density is only implicitly specified by the complete data model, $f_\theta(y) = \int f_\theta(x, y) \, dx$, and is often not available in closed form. This is what makes likelihood inference for missing data difficult.

Many Monte Carlo methods for approximating the observed data likelihood in a missing data model have been proposed. In these, missing data are simulated by either ordinary Monte Carlo [17, 24] or Markov chain Monte Carlo (MCMC) [10, 12, 18, 28, 29]. To get a good approximation of the likelihood over a large region, umbrella sampling [13, 30] may be necessary. There are also many Monte Carlo methods for maximum likelihood without approximating the likelihood: stochastic approximation [23, 38], Monte Carlo EM [14, 34], and Monte Carlo Newton–Raphson [25]. There are also non-Monte Carlo methods for maximum likelihood without approximating the likelihood: EM [8] and analytic approximation [6]. There are so many methods because each has its strength and weakness. In theory, Monte Carlo methods work for complicated problems but require very careful calibration, whereas non-Monte Carlo methods are relatively easier to implement but only apply to simple cases. All are useful for some, but not all, problems.

This article is concerned with a Monte Carlo approximation of the observed data likelihood and asymptotic properties of the maximizer of our Monte Carlo likelihood. Our method uses simulated missing data that are independent and identically distributed (i.i.d.) and independent of the observed data. It approximates the likelihood over the entire parameter space. We give Logit–Normal GLMM examples to illustrate the case when our asymptotic normality holds and to show the value of our Monte Carlo likelihood approximation even when asymptotic normality does not hold.

Let the observed data $Y_1, \ldots, Y_n$ be i.i.d. from a density $g$, which is not assumed to be some $f_\theta$. We do not assume the model is correctly specified, since this increase of generality makes the theory no more difficult. The MLE $\hat{\theta}_n$ is a maximizer of the log likelihood

$$l_n(\theta) = \sum_{j=1}^{n} \log f_\theta(Y_j). \tag{1}$$

In our method, we generate an i.i.d. Monte Carlo sample $X_1, \ldots, X_m$, independent of $Y_1, \ldots, Y_n$, from an importance sampling density $h$ and approximate $f_\theta(y)$ by

$$f_{\theta,m}(y) = \frac{1}{m} \sum_{i=1}^{m} \frac{f_\theta(X_i, y)}{h(X_i)}. \tag{2}$$

This makes heuristic sense because

$$f_{\theta,m}(y) \xrightarrow{\text{a.s.}}_m E_h\left\{\frac{f_\theta(X,y)}{h(X)}\right\} = f_\theta(y) \qquad \text{for each } y$$

by the strong law of large numbers. (The subscript $m$ on the arrow means as $m$ goes to infinity. Similarly, a subscript $m, n$ means as both $m$ and $n$ go to infinity.) Our estimate of $\hat{\theta}_n$ is the maximizer $\hat{\theta}_{m,n}$ of our Monte Carlo log likelihood

$$l_{m,n}(\theta) = \sum_{j=1}^{n} \log f_{\theta,m}(Y_j), \tag{3}$$

an approximation to $l_n(\theta)$ with $f_{\theta,m}$ replacing $f_\theta$. We call $\hat{\theta}_{m,n}$ the Monte Carlo MLE (MCMLE).

Under the conditions of Theorem 2.3, the MCMLE

$$\hat{\theta}_{m,n} \approx \mathcal{N}\left(\theta^*, \frac{J^{-1}VJ^{-1}}{n} + \frac{J^{-1}WJ^{-1}}{m}\right), \tag{4}$$

for sufficiently large Monte Carlo sample size $m$ and observed data sample size $n$, where $\theta^*$ is the minimizer of the Kullback–Leibler information

$$K(\theta) = E_g \log \frac{g(Y)}{f_\theta(Y)}, \tag{5}$$

$J$ is minus the expectation of the second derivative of the log likelihood, $V$ is the variance of the first derivative of the log likelihood (score), and $W$ is the variance of the deviation of the score from its Monte Carlo approximation [given by (7) below]. Under certain regularity conditions [15, 35],

$$\hat{\theta}_n \approx \mathcal{N}\left(\theta^*, \frac{J^{-1}VJ^{-1}}{n}\right). \tag{6}$$

We see that $\hat{\theta}_{m,n}$ has nearly the same distribution when the Monte Carlo sample size $m$ is very large. If the model is correctly specified, that is, $g = f_{\theta_0}$, then $\theta^* = \theta_0$ and $J = V$, either of which is called Fisher information, and (6) becomes

$$\hat{\theta}_n \approx \mathcal{N}\left(\theta^*, \frac{J^{-1}}{n}\right),$$

the familiar formula due to Fisher and Cramér. This replacement of $J^{-1}$ by the so-called "sandwich" $J^{-1}VJ^{-1}$ is the only complication arising from model misspecification [19].

The first term of the asymptotic variance in (4) is what would be the asymptotic variance if we could use the exact likelihood rather than Monte Carlo. Hence it is the same as the asymptotic variance in (6). The second term is additional variance due to Monte Carlo. Increasing the Monte Carlo sample size $m$ can make the second term as small we please so that the MCMLE $\hat{\theta}_{m,n}$ is almost as good as the MLE $\hat{\theta}_n$. In (4), $W$ is the only term related to the importance sampling density $h$ that generates the Monte Carlo sample. Choosing an $h$ that makes $W$ smaller makes $\hat{\theta}_{m,n}$ more accurate.

The asymptotic distribution of $\hat{\theta}_{m,n}$ in (4) is a convolution of two independent normal distributions. The proof of this is not simple, however, for three reasons. First, the finite sample terms from which these arise [the two terms on the right-hand side in (9) below] are dependent. Second, one of these is itself a sum of dependent terms, because each term in (3) uses the same $X$'s. Third, our two sample sizes $m$ and $n$ go to infinity simultaneously, and we must show that the result does not depend on the way in which $m$ and $n$ go to infinity.

**2. Asymptotics of $\hat{\theta}_{m,n}$.** In this section, we state theorems about strong consistency and asymptotic normality of the MCMLE $\hat{\theta}_{m,n}$. Proofs are in the Appendix.

We use empirical process notation throughout. We let $P$ denote the probability measure induced by the importance sampling density $h$, and we let $\mathbb{P}_m$ denote the empirical measure induced by $X_1,\ldots,X_m$ (that are i.i.d. from $P$). Similarly, we let $Q$ denote the probability measure induced by the true density $g$ and $\mathbb{Q}_n$ denote the empirical measure induced by $Y_1,\ldots,Y_n$ (that are i.i.d. from $Q$). Given a measurable function $f:\mathcal{X} \mapsto \mathbb{R}$, we write $\mathbb{P}_m f(X)$ for the expectation of $f$ under $\mathbb{P}_m$ and $Pf(X)$ for the expectation under $P$. Similarly we use $\mathbb{Q}_n f(Y)$ and $Qf(Y)$. Note that $\mathbb{P}_m f(X) = \frac{1}{m}\sum_{i=1}^m f(X_i)$ is just another notation for a particular sample mean.

The Kullback–Leibler information in (5) is written as $K(\theta) = Q\log[g(Y)/f_\theta(Y)]$, its empirical version as $K_n(\theta) = \mathbb{Q}_n \log[g(Y)/f_\theta(Y)]$ and our approximation to $K_n(\theta)$ as

$$K_{m,n}(\theta) = \mathbb{Q}_n \log[g(Y)/f_{\theta,m}(Y)]$$

with $f_{\theta,m}(y) = \mathbb{P}_m f_\theta(X,y)/h(X)$. Then $K_n(\theta) = \mathbb{Q}_n \log g(Y) - l_n(\theta)/n$ and $K_{m,n}(\theta) = \mathbb{Q}_n \log g(Y) - l_{m,n}(\theta)/n$. Hence the MLE $\hat{\theta}_n$ is the minimizer of $K_n$ and the MCMLE $\hat{\theta}_{m,n}$ is the minimizer of $K_{m,n}$. By Jensen's inequality $K(\theta) \geq 0$. This allows $K(\theta) = \infty$ for some $\theta$, but we assume $K(\theta^*)$ is finite. [This excludes only the uninteresting case of the function $\theta \mapsto K(\theta)$ being identically $\infty$.]

2.1. *Epi-convergence of $K_{m,n}$.* To get the convergence of $\hat{\theta}_{m,n}$ to $\theta^*$ we use *epi-convergence* of the function $K_{m,n}$ to the function $K$. Epi-convergence is a "one-sided" uniform convergence that was first introduced by Wijsman [36, 37], developed in optimization theory [2, 3, 26] and used in statistics [11, 12]. It is weaker than uniform convergence yet insures the convergence of minimizers as the following proposition due to Attouch [2], Theorem 1.10, describes.

PROPOSITION 2.1. *Let $X$ be a general topological space, $\{f_n\}$ a sequence of functions from $X$ to $\overline{\mathbb{R}}$ that epi-converges to $f$, and $\{x_n\}$ a sequence of points in $X$ satisfying $f_n(x_n) \leq \inf f_n + \varepsilon_n$ with $\varepsilon_n \downarrow 0$. Then for every converging subsequence $x_{n_k} \to x_0$*

$$f(x_0) = \inf f = \lim_k f_{n_k}(x_{n_k}).$$

If $f$ has a unique minimizer $x$, then $x$ is the only cluster point of the sequence $\{x_n\}$. Otherwise, there may be many cluster points, but all must minimize $f$. There may not be any convergent subsequence. If the sequence $\{x_n\}$ is in a compact set and $X$ is sequentially compact, however, there is always a convergent subsequence.

THEOREM 2.2. *Let $\{f_\theta(x,y): \theta \in \Theta\}$, where $\Theta \subset \mathbb{R}^d$, be a family of densities with respect to a $\sigma$-finite measure $\mu \times \nu$ on $\mathcal{X} \times \mathcal{Y}$, let $X_1, X_2, \ldots$ be i.i.d. from a probability distribution $P$ that has a density $h$ with respect to $\mu$, and let $Y_1, Y_2, \ldots$ be i.i.d. from a probability distribution $Q$ that has a density $g$ with respect to $\nu$. Suppose:*

(1) *$\Theta$ is a second countable topological space;*
(2) *for each $(x,y)$, the function $\theta \mapsto f_\theta(x,y)$ is upper semicontinuous on $\Theta$;*
(3) *for each $\theta$, there exists a neighborhood $B_\theta$ of $\theta$ such that*

$$Q \log\left[P \sup_{\phi \in B_\theta} f_\phi(X,Y)/h(X)g(Y)\right] < \infty;$$

(4) *for each $\theta$, there exists a neighborhood $C_\theta$ of $\theta$ such that for any subset $B$ of $C_\theta$, the family of functions $\{\sup_{\phi \in B} f_\phi(\cdot, y)/h(\cdot)g(y) : y \in \mathcal{Y}\}$ is $P$-Glivenko–Cantelli;*
(5) *for each $\theta$, the family of functions $\{f_\theta(\cdot|y)/h(\cdot) : y \in \mathcal{Y}\}$ is $P$-Glivenko–Cantelli.*

*Then $K_{m,n}$ epi-converges to $K$ with probability one.*

Glivenko–Cantelli means a family of functions for which the uniform strong law of large numbers holds ([32], page 81). Conditions (1) through (3) are similar to those of Theorem 2 in [12]. Also they are vaguely similar to those in [33], which imply epi-convergence of $K_n$ to $K$ (when there are no missing data and no Monte Carlo).

2.2. *Asymptotic normality of $\hat{\theta}_{m,n}$.* The following theorem assumes that the local minimizer $\theta^*$ of $K$ is an interior point of $\Theta$ and that $K$ is differentiable. Hence $\nabla K(\theta^*) = 0$, where $\nabla$ means differentiation with respect to $\theta$.

THEOREM 2.3. *Let $\{f_\theta(x,y): \theta \in \Theta\}$, where $\Theta \subset \mathbb{R}^d$, be a family of densities with respect to a $\sigma$-finite measure $\mu \times \nu$ on $\mathcal{X} \times \mathcal{Y}$, let $X_1, X_2, \ldots$ be i.i.d. from a probability distribution $P$ that has a density $h$ with respect to $\mu$, and let $Y_1, Y_2, \ldots$ be i.i.d. from a probability distribution $Q$ that has a density $g$ with respect to $\nu$. Suppose:*

(1) *second partial derivatives of $f_\theta(x,y)$ with respect to $\theta$ exist and are continuous on $\Theta$ for all $x$ and $y$, and may be passed under the integral sign in $\int f_\theta(x|y)\, d\mu(x)$;*
(2) *$\mathcal{Y}$ is a separable metric space and $y \mapsto \nabla f_{\theta^*}(x|y)$ is continuous for each $x$;*

(3) *there is an interior point $\theta^*$ of $\Theta$ such that $Q\nabla \log f_{\theta^*}(Y) = 0$, $V = \mathrm{var}_Q \nabla \log f_{\theta^*}(Y)$ is finite and $J = -Q\nabla^2 \log f_{\theta^*}(Y)$ is finite and nonsingular;*

(4) *there exists a $\rho > 0$ such that $S_\rho = \{\theta : |\theta - \theta^*| \leq \rho\}$ is contained in $\Theta$ and $\mathcal{F}_1 = \{\nabla^2 f_\theta(\cdot) : \theta \in S_\rho\}$ is $Q$-Glivenko–Cantelli;*

(5) *$\mathcal{F}_2 = \{f_{\theta^*}(\cdot|y)/h(\cdot) : y \in \mathcal{Y}\}$ is $P$-Glivenko–Cantelli;*

(6) *$\mathcal{F}_3 = \{\nabla f_{\theta^*}(\cdot|y)/h(\cdot) : y \in \mathcal{Y}\}$ is $P$-Donsker and its envelope function $F$ has a finite second moment;*

(7) *$\mathcal{F}_4 = \{\nabla^2 f_\theta(\cdot|y)/h(\cdot) : y \in \mathcal{Y}, \theta \in S_\rho\}$ is $P$-Glivenko–Cantelli;*

(8) *there is a sequence $\hat{\theta}_{m,n}$ which converges to $\theta^*$ in probability such that*

$$\sqrt{\min(m,n)}\nabla K_{m,n}(\hat{\theta}_{m,n}) \xrightarrow{P}_{m,n} 0.$$

*Then*

(7) $$W = \mathrm{var}_P Q\nabla f_{\theta^*}(X|Y)/h(X)$$

*is finite and*

(8) $$\left(\frac{V}{n} + \frac{W}{m}\right)^{-1/2} J(\hat{\theta}_{m,n} - \theta^*) \xrightarrow{\mathcal{L}}_{m,n} \mathcal{N}(0, I).$$

Donsker means a family of functions for which the uniform central limit theorem holds ([32], page 81). Note $\mathcal{F}_3$ is a family of vector-valued functions and $\mathcal{F}_1$ and $\mathcal{F}_4$ are families of matrix-valued functions. Such families are Glivenko–Cantelli or Donsker if each component is ([31], page 270). Conditions (1), (3), (4) and (8) are similar to the usual regularity conditions for asymptotic normality of the MLE, which can be found, for example, in [9], Chapter 18. For a correctly specified model, differentiability under the integral sign in $1 = \iint f_\theta(x,y)\,d\mu(x)\,d\nu(y)$ implies conditions (1) and (3). Condition (4) holds if functions in $\mathcal{F}_1$ are dominated by a $L^1(Q)$ function because $S_\rho$ is compact ([9], Theorem 16(a)).

Under smoothness conditions imposed in this theorem, the asymptotics of $\hat{\theta}_{m,n}$ arises from the asymptotics of

(9) $$\nabla K_{m,n}(\theta^*) = -\mathbb{Q}_n \nabla \log f_{\theta^*}(Y) - \mathbb{Q}_n \nabla \log \mathbb{P}_m f_{\theta^*}(X|Y)/h(X).$$

The two terms on the right-hand side are dependent and the summands in the second term are dependent, which indicates the complexity of this problem and why the usual asymptotic arguments do not work here. The asymptotics for the first term follow from the central limit theorem. The asymptotics for the second term (Lemma A.4) go as shown below:

$$\sqrt{m}\mathbb{Q}_n \nabla \log \mathbb{P}_m f_{\theta^*}(X|Y)/h(X) \xrightarrow{m} \mathbb{Q}_n \mathbb{G}_P \nabla f_{\theta^*}(X|Y)/h(X)$$
$$\searrow_{m,n} \qquad \downarrow_n$$
$$Q\mathbb{G}_P \nabla f_{\theta^*}(X|Y)/h(X)$$

We first let $m \to \infty$ then $n \to \infty$. A uniformity argument then makes the result the same when $m$ and $n$ go to infinity simultaneously. The $\xrightarrow{m}$ part is weak convergence of the empirical process $\sqrt{m}(\mathbb{P}_m - P)$ to a tight Gaussian process $\mathbb{G}_P$. The $\xrightarrow{n}$ part is the law of large numbers. Integration over sample paths of $\mathbb{G}_P$ gives the distribution of the limit. The asymptotic independence between the two terms in (9) comes from the fact that the law of large numbers eliminates the randomness coming from $\mathbb{Q}_n$.

2.3. *Plug-in estimates for $J$, $V$ and $W$.* We can construct a confidence region for $\theta^*$ using (4) or (8). If we can evaluate the integrals defining $J$, $V$ and $W$, then we may use those integrals with $\hat{\theta}_{m,n}$ plugged in for $\theta^*$ to estimate them, assuming enough continuity. Often we cannot evaluate the integrals or do not know $g$. Then we use their sample versions,

$$\widehat{J}_{m,n} = -\frac{1}{n} \sum_{j=1}^{n} \nabla^2 \log f_{\hat{\theta}_{m,n}}(Y_j),$$

(10) $$\widehat{V}_{m,n} = \frac{1}{n} \sum_{j=1}^{n} \nabla \log f_{\hat{\theta}_{m,n}}(Y_j) \nabla \log f_{\hat{\theta}_{m,n}}(Y_j)^T,$$

$$\widehat{W}_{m,n} = \frac{1}{m} \sum_{i=1}^{m} \widehat{S}_i \widehat{S}_i^T,$$

where

(11) $$\widehat{S}_i = \frac{1}{n} \sum_{j=1}^{n} \nabla f_{\hat{\theta}_{m,n}}(X_i|Y_j)/h(X_i).$$

Often these cannot be used as shown because $f_\theta(y)$ and $f_\theta(x|y)$ are not available in closed form. Then we replace $f_\theta(y)$ by $f_{\theta,m}(y)$ defined in (2) and $f_\theta(x|y)$ by $f_\theta(x,y)/f_{\theta,m}(y)$. The resulting variance estimate $\widehat{J}_{m,n}^{-1}(\widehat{V}_{m,n}/n + \widehat{W}_{m,n}/m)\widehat{J}_{m,n}^{-1}$ is the sandwich estimator.

2.4. *An alternative Monte Carlo scheme.* Each term in (3) uses the same $X$'s. An alternative is to use each $X$ once, generating a new sample for each term in (3). Then the resulting estimate has the same asymptotic variance as in (4) or (8) except that $W$ is replaced by $\widetilde{W} = Q \operatorname{var}_h \nabla f_{\theta^*}(X|Y)/h(X)$. By Jensen's inequality, $\widetilde{W} \geq W$. Thus using the $X$'s $n$ times makes $\hat{\theta}_{m,n}$ more accurate.

**3. Logit–Normal GLMM examples.** The Logit–Normal GLMM refers to Bernoulli regression with normal random effects. It has a linear predictor of the form

(12) $$\eta = X\beta + Zb,$$

where $X$ and $Z$ are known design matrices, and $\beta$ and $b$ are unknown vectors (fixed effects and random effects, resp.). The observed data consist of $n$ i.i.d. responses, one for each individual, and the missing data consist of $n$ i.i.d. random effects vectors, one for each individual with $b \sim \mathcal{N}(0, \Sigma)$ (we denote the missing data by $b$, instead of $x$, to avoid confusion with $X$). The observed data for one individual is a vector whose components are independent Bernoulli given $b$, with success probability vector having components $\text{logit}^{-1}(\eta_k) = 1/(1 + \exp(-\eta_k))$. The unknown parameters to be estimated are $\beta$ and the parameters determining the variance $\Sigma$ of random effects, which typically has simple structure and involves only a few parameters.

We reparametrized (12) as

$$\eta = X\beta + Z\Delta b, \tag{13}$$

where $\Delta$ is a diagonal matrix whose diagonal is a vector of unknown parameters (square roots of variance components) and $b$ is a standard normal random vector (whose distribution contains no unknown parameters). All of the unknown parameters are in $\beta$ and $\Delta$ in the linear predictor (13). This representation is flexible enough to include the examples in this article. We used the standard normal density (which is the true density of $b$) as our importance sampling density. This makes sense because of our reparametization to make the density of $b$ not depend on the parameters. This is not a general recommendation of the normal density.

We wrote an R package `bernor` that implements the methods of this article for the Logit–Normal GLMM (available at www.stat.umn.edu/geyer/bernor). The web page also contains detailed verification of the conditions of our theorems for the model and detailed descriptions of its applications to our examples.

3.1. *Conditions of the theorems.* The Logit–Normal GLMM with our importance sampling density satisfies the conditions of both theorems. Verifying the conditions is straightforward because of two properties. First, the sample space $\mathcal{Y}$ is finite. Thus verifying Glivenko–Cantelli in conditions (4) and (5) of Theorem 2.2 and condition (6) of Theorem 2.3 reduces to just verifying that functions are $L^1(P)$, and verifying Donsker in condition (5) of Theorem 2.3 reduces to just verifying that functions are $L^2(P)$. Also verifying Glivenko–Cantelli in condition (7) of Theorem 2.3 reduces to just verifying that for each $y$, the class $\{\nabla^2 f_\theta(\cdot|y)/h(\cdot) : \theta \in S_\rho\}$ is $P$-Glivenko–Cantelli. This can be verified like condition (4) of Theorem 2.3 as discussed after the theorem. Second, our importance sampling density $h$ is the marginal density of the missing data and this implies $f_\theta(b,y)/h(b) = f_\theta(y|b)$, which makes it easy to verify that functions are $L^1$ or $L^2$. Differentiability under the integral sign twice follows from $h$ having two moments.

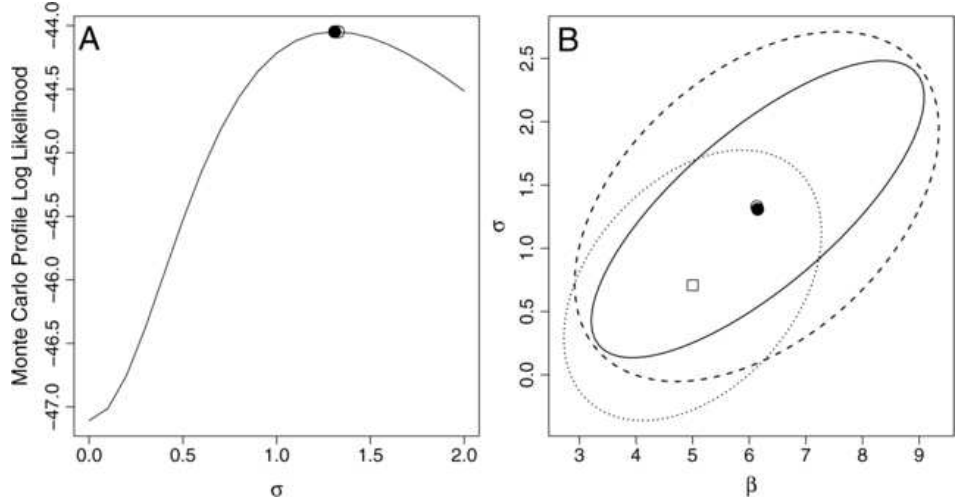

FIG. 1. *Monte Carlo profile log likelihood* (A) *and nominal* 95% *confidence ellipses* (B) *for the Booth and Hobert data using* $m = 10^4$. *Solid dot and solid line are the MCMLE and confidence ellipse using plug-in estimates of* $J$, $V$ *and* $W$ *at the MCMLE. Hollow dot and dashed line are the MLE and confidence ellipse using Fisher information and exact* $W$ *at the MLE. Square and dotted line are the "simulation truth" parameter value and confidence ellipse using Fisher information and exact* $W$ *at the simulation truth. The last two assume* $V = J$.

3.2. *Data from McCulloch's model.* We use a data set given by Booth and Hobert [5], Table 2, that was simulated using a model from [22]. This model corresponds to a Logit–Normal GLMM with one-dimensional $\beta$ and $b$ in (12), and its log likelihood can be calculated exactly by numerical integration. The observed data consist of ten i.i.d. vectors of length 15. The parameters that generated the data are $\beta = 5$ and $\sigma = \sqrt{1/2}$.

Using a Monte Carlo sample of size $10^4$, we approximated the observed data log likelihood and obtained the MCMLE. The Monte Carlo profile log likelihood for $\sigma$ (Figure 1A) indicates that the log likelihood is well behaved, quadratic around the MLE, and that our MCMLE ($\hat{\beta}_{m,n} = 6.15, \hat{\sigma}_{m,n} = 1.31$) is very close to the MLE ($\hat{\beta}_n = 6.13, \hat{\sigma}_n = 1.33$). Using plug-in estimates given by (10), we also obtained a nominal 95% confidence ellipse for the true parameter (the solid ellipse in Figure 1B). For comparison, we obtained two other confidence ellipses using the theoretical expected Fisher information and $W$ at the MLE (the dashed ellipse in Figure 1B) and also at the true parameter (the dotted ellipse in Figure 1B). Both these exact evaluations took 13 hours, whereas our plug-in estimates took two and-a-half minutes. Our MCMLE and the MLE are not close to the truth, and these ellipses are different, indicating that an observed data sample size $n = 10$ is too small to apply asymptotics. But our MCMLE is close to the MLE, indicating that our Monte Carlo sample size $m = 10^4$ is good enough for estimating the MLE for the observed data.

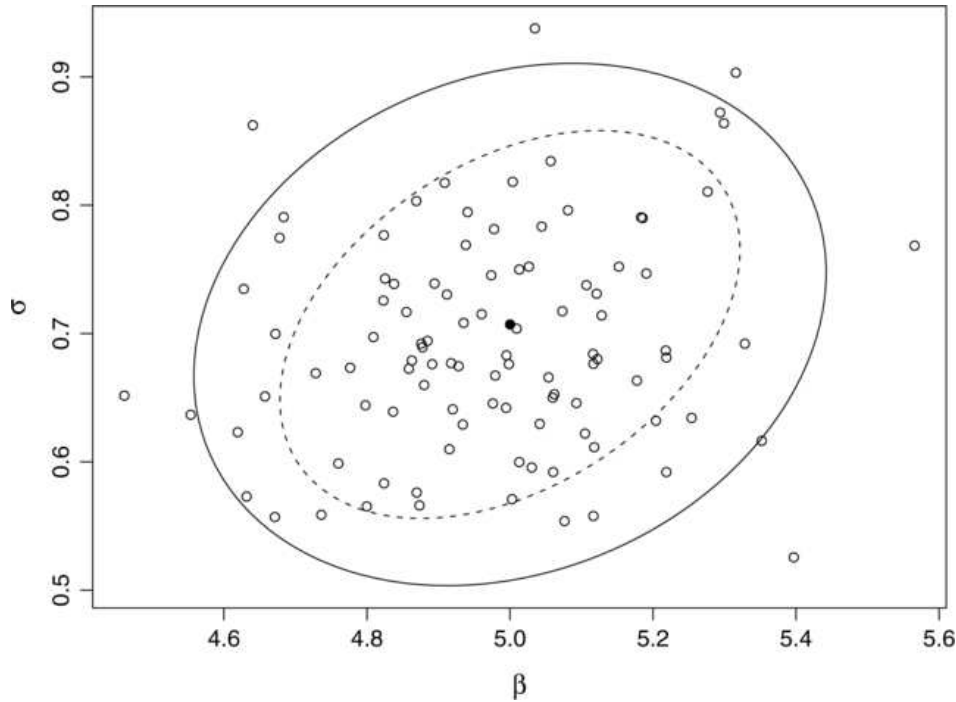

Fig. 2. *Sampling distribution of the MCMLE. Hollow dots are the MCMLE's for* 100 *simulated data sets, using sample sizes* $n = 500$ *and* $m = 100$. *The solid dot is the "simulation truth" parameter value. The solid curve is the asymptotic* 95% *coverage ellipse. The dashed curve is what the* 95% *coverage ellipse would be if* $m$ *were infinity.*

3.3. *Simulation for McCulloch's model.* To demonstrate our asymptotic theory, we did a simulation study using the same model with sample sizes $n = 500$ and $m = 100$. [We chose these sample sizes so that the two terms that make up the variance in (4) have roughly the same size.] Figure 2 gives the scatter plot of 100 MCMLE's. The solid ellipse is an asymptotic 95% coverage ellipse using the theoretical expected Fisher information and $W$. The dashed ellipse is what we would have if we had very large Monte Carlo sample size $m$, leaving $n$ the same. The solid ellipse contains 92 out of 100 points, thus asymptotics appear to work well at these sample sizes.

3.4. *The influenza data.* Table 1 in [7] shows data collected from 263 individuals about four influenza outbreaks from 1977 to 1981 in Michigan. Thus the observed data consist of 263 i.i.d. vectors of length four. Coull and Agresti [7] used a Logit–Normal GLMM with four-dimensional $\beta$ and $b$ in

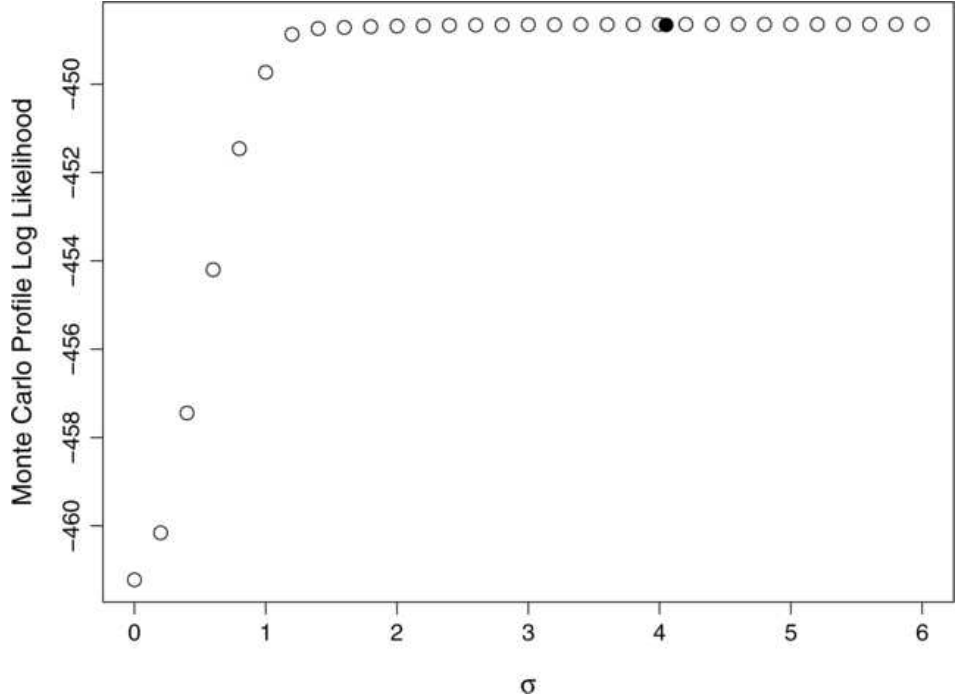

FIG. 3. *Monte Carlo profile log likelihood using $m = 10^6$. For each $\sigma$, other parameters are maximized. The solid dot is the MLE reported by Coull and Agresti* [7]. *Leftmost point ($\sigma = 0$) corresponds to the MLE for the model without random effects.*

(12) and $b$ having variance matrix

$$\sigma^2 \begin{bmatrix} 1 & \rho_1 & \rho_1 & \rho_2 \\ \rho_1 & 1 & \rho_1 & \rho_2 \\ \rho_1 & \rho_1 & 1 & \rho_2 \\ \rho_2 & \rho_2 & \rho_2 & 1 \end{bmatrix}$$

and reported the MLE as $\hat{\beta} = (-4.0, -4.4, -4.7, -4.5)$, $\hat{\sigma} = 4.05$, $\hat{\rho_1} = 0.43$ and $\hat{\rho_2} = -0.25$. Our reparametrization (13) that corresponds to this model has four-dimensional identity matrix $X$, four-dimensional $\beta$,

$$Z = \begin{pmatrix} 1 & 1 & 1 & 0 & 0 & 0 \\ 1 & 1 & 0 & 1 & 0 & 0 \\ 1 & 1 & 0 & 0 & 1 & 0 \\ 1 & -1 & 0 & 0 & 0 & 1 \end{pmatrix},$$

six-dimensional diagonal matrix $\Delta$ with diagonal elements $\delta_1, \delta_2, \delta_3, \delta_3, \delta_3, \delta_3$, and six-dimensional $b$.

Using a Monte Carlo sample of size $10^6$, we approximated the observed data log likelihood and found a ridge in the log likelihood surface (Figure 3).

TABLE 1
*Parameter values along the ridge of the Monte Carlo log likelihood for the Influenza Data, using Monte Carlo sample size $m = 10^6$. The MLE from [7] is provided in the last row for comparison*

| $\sigma$ | $\rho_1$ | $\rho_2$ | $\beta$ | | | | MC log likelihood |
|------|------|--------|-------|-------|-------|-------|-----------|
| 1.60 | 0.79 | $-0.47$ | $-2.1$ | $-2.3$ | $-2.5$ | $-2.4$ | $-448.717$ |
| 2.00 | 0.64 | $-0.38$ | $-2.4$ | $-2.6$ | $-2.8$ | $-2.7$ | $-448.682$ |
| 3.00 | 0.48 | $-0.28$ | $-3.2$ | $-3.5$ | $-3.7$ | $-3.6$ | $-448.646$ |
| 4.00 | 0.43 | $-0.25$ | $-4.0$ | $-4.4$ | $-4.6$ | $-4.5$ | $-448.635$ |
| 5.00 | 0.40 | $-0.23$ | $-4.8$ | $-5.3$ | $-5.6$ | $-5.5$ | $-448.631$ |
| 6.00 | 0.39 | $-0.22$ | $-5.7$ | $-6.2$ | $-6.6$ | $-6.4$ | $-448.629$ |
| 4.05 | 0.43 | $-0.25$ | $-4.0$ | $-4.4$ | $-4.7$ | $-4.5$ | $-448.646$ |

(Monte Carlo sample size $10^7$ gave results identical to three decimal places.) The log likelihood is strongly curved in directions orthogonal to the ridge but hardly changes along the ridge. Fisher information is nearly singular because of this ridge. Parameter values along the ridge (Table 1) vary over a large range, and the bigger $\sigma$ is the more extreme the components of $\beta$ are. This is a surprise because sample size 263 is usually large enough for making inference about seven parameters. Even though the model is identifiable, it is not clear that asymptotics would hold for any sample size. Hence some penalized likelihood method should probably be used.

3.5. *The salamander data.* We use the data in [21], Section 14.5, that were obtained from a salamander mating experiment and have been analyzed many times (see [5], for one analysis and citations of others). This example has been considered difficult to analyze because its likelihood involves a 20-dimensional integral. We use "Model A" of Karim and Zeger [16], which corresponds to a Logit–Normal GLMM with four-dimensional $\beta$ and 20-dimensional $b$ in (12) with two parameters determining the variance of $b$. The observed data consist of three i.i.d. vectors of length 120. The MLE given by Booth and Hobert [5] is $\hat{\beta} = (1.03, 0.32, -1.95, 0.99)$ and $\hat{\sigma} = (1.18, 1.12)$. Based on Monte Carlo sample size $10^7$, our MCMLE was $\hat{\beta}_{m,n} = (1.00, 0.53, -1.78, 1.27)$ and $\hat{\sigma}_{m,n} = (1.10, 1.17)$, and the standard errors were (0.35, 0.33, 0.36, 0.53) for $\hat{\beta}_{m,n}$ and (0.20, 0.28) for $\hat{\sigma}_{m,n}$. Our method did not work well for these data, and these standard errors give a clear indication of the accuracy of our MCMLE.

**4. Discussion.** We have described a Monte Carlo method to approximate the observed data likelihood and the MLE when there are missing data and the observed data likelihood is not available in closed form. The MLE converges to the minimizer $\theta^*$ of the Kullback–Leibler information, which is the true parameter value when the model is correctly specified. We have proved that our MCMLE is a consistent and asymptotically normal estimate of $\theta^*$ as both Monte Carlo and observed data sample sizes go to infinity simultaneously. Plug-in estimates of the asymptotic variance are provided in (10) for constructing confidence regions for $\theta^*$.

We have presented the theory so that it can be used for studying model misspecification in missing data models. In practice, a statistical model $f_\theta$ is often chosen only for mathematical convenience and may contain simplistic and unrealistic assumptions. However, it is usually possible to simulate i.i.d. data $Y$'s from a more realistic model $g$. The theory applies whether the $Y$'s are a Monte Carlo sample or real data. In either case we can estimate $\theta^*$ using $\hat{\theta}_{m,n}$ and know what accuracy we have. By comparing $f_{\hat{\theta}_{m,n}}$ (an estimate of $f_{\theta^*}$, the "best" approximation to $g$ in the model) with $g$, we can assess model validity as whether the particular model is reasonable for approximating the truth or how its simplifying assumptions influence scientific conclusions.

Our applications to the Logit–Normal GLMM examples illustrate advantages and disadvantages of our method. First, our method uses ordinary (independent sample) Monte Carlo, thus is simpler to implement and easier to understand than MCMC. Second, it always provides accurate standard errors, and they give a clear indication of when the method works and how well. MCMC methods require more careful tuning and do not provide analogous standard errors. MCMC diagnostics are widely used but give no guarantees, and convergence proofs are very difficult except for simple applications and are not widely used. Third, our method approximates the likelihood over the entire parameter space. We have seen the advantage of such likelihood evaluation for the influenza data in Section 3.4. One can assess whether the likelihood is well behaved so that appropriate inference can be based on the MLE. The only disadvantage of our method arises from its simple Monte Carlo scheme. It does not work well with high-dimensional missing data as in the salamander data in Section 3.5.

Our method is based on sampling from an importance sampling density $h$. In theory, we want the optimal $h$ that makes $W$ as small as possible so that the MCMLE is as accurate as possible. The form of $W$ in (7) says that we want $h(x)$ to be high where $Q\nabla f_{\theta^*}(x|Y)$ is high. In very simple situations we can find such $h$ (Sung [27] finds the optimal $h$ for a normal mixture model). In complicated situations, just as in ordinary importance sampling, one cannot calculate the optimal $h$ and must proceed by trial and error.

Asymptotic theory analogous to ours does not exist for MCMC. It involves three quantities: the MCMLE $\hat{\theta}_{m,n}$ is a function of both simulated missing data and observed data, the MLE $\hat{\theta}_n$ (which cannot be calculated exactly) is a function of observed data only, and $\theta^*$ is the true parameter value. Geyer [12] provides asymptotic theory for $\hat{\theta}_{m,n} - \hat{\theta}_n$ conditional on observed data, accounting for only Monte Carlo variability, not sampling variability. Classical theory of maximum likelihood provides asymptotic theory for $\hat{\theta}_n - \theta^*$, accounting for sampling variability. As we have seen in this article, it is not easy to combine these two sources of variability, and this has not been tried for MCMC. Our method could be extended so that the importance sampling density can depend on observed data, which is usually done in MCMC. We suppose the theory for that would be considerably more complicated than what we have presented here and would be even more complicated for MCMC.

Even though our original motivation was theoretical, our method does work in practical examples. The `bernor` package can be used for analysis of Logit–Normal GLMM. Our method is applicable to other missing data models.

## APPENDIX

**A.1. Proof of Theorem 2.2.** Let $(m_k, n_k)$ be a subsequence. We need to show

$$K \leq \text{e-}\liminf_k K_{m_k,n_k} \leq \text{e-}\limsup_k K_{m_k,n_k} \leq K,$$

which is equivalent to

$$(14) \qquad K(\theta) \leq \sup_{B \in \mathcal{N}(\theta)} \liminf_{k \to \infty} \inf_{\phi \in B} K_{m_k,n_k}(\phi),$$

$$(15) \qquad K(\theta) \geq \sup_{B \in \mathcal{N}(\theta)} \limsup_{k \to \infty} \inf_{\phi \in B} K_{m_k,n_k}(\phi),$$

where $\mathcal{N}(\theta)$ is the set of neighborhoods of $\theta$. By condition (1) there is a countable basis $\mathcal{B} = \{B_1, B_2, \ldots\}$ for the topology of $\Theta$. Choose a countable dense subset $\Theta_c = \{\theta_1, \theta_2, \ldots\}$ by choosing $\theta_n \in B_n$ to satisfy $K(\theta_n) \leq \inf_{\phi \in B_n} K(\phi) + 1/n$. Let $\mathcal{N}_c(\theta) = \{B \in \mathcal{B} \cap \mathcal{N}(\theta) : B \subset B_\theta \cap C_\theta\}$ where $B_\theta$ is given by condition (3) and $C_\theta$ is given by condition (4). Suprema over $\mathcal{N}(\theta)$ in (14) and (15) can be replaced by suprema over the countable set $\mathcal{N}_c(\theta)$.

By Lemma A.1 below

$$(16) \qquad \limsup_{k \to \infty} K_{m_k,n_k}(\theta) \leq K(\theta)$$

for each $\theta$ with probability one, and by Lemma A.2 below

$$(17) \qquad \liminf_{k \to \infty} \inf_{\phi \in B} K_{m_k,n_k}(\phi) \geq -Q \log P \sup_{\phi \in B} \frac{f_\phi(X,Y)}{h(X)g(Y)}$$

for each $B \in \mathcal{N}_c(\theta)$ with probability one. Since $\Theta_c$ and $\bigcup_{\theta \in \Theta} \mathcal{N}_c(\theta)$ are countable and since a countable union of null sets is a null set, we have (16) and

(17) simultaneously on $\Theta_c$ and $\bigcup_{\theta\in\Theta}\mathcal{N}_c(\theta)$ with probability one. If $B\in\mathcal{B}$ and $\theta\in B\cap\Theta_c$, then by (16)

$$K(\theta)\geq\limsup_k K_{m_k,n_k}(\theta)\geq\limsup_{k\to\infty}\inf_{\phi\in B}K_{m_k,n_k}(\phi).$$

Hence

$$\sup_{B\in\mathcal{N}_c(\theta)}\inf_{\phi\in B\cap\Theta_c}K(\phi)\geq\sup_{B\in\mathcal{N}_c(\theta)}\limsup_{k\to\infty}\inf_{\phi\in B}K_{m_k,n_k}(\phi).$$

The term on the left-hand side is $K(\theta)$ by lower semicontinuity of $K$ (Lemma A.3 below) and by the construction of $\Theta_c$. This proves (15). We also have

$$\sup_{B\in\mathcal{N}_c(\theta)}\liminf_{k\to\infty}\inf_{\phi\in B}K_{m_k,n_k}(\phi)\geq\sup_{B\in\mathcal{N}_c(\theta)}-Q\log P\sup_{\phi\in B}\frac{f_\phi(X,Y)}{h(X)g(Y)}$$

$$=-Q\log P\inf_{B\in\mathcal{N}_c(\theta)}\sup_{\phi\in B}\frac{f_\phi(X,Y)}{h(X)g(Y)}$$

$$=-Q\log P\frac{f_\theta(X,Y)}{h(X)g(Y)}=K(\theta),$$

where the inequality follows from (17), the first equality from the monotone convergence theorem and the second equality from condition (2). This proves (14).

LEMMA A.1. *Under condition* (5) *of Theorem* 2.2, $K_{m,n}(\theta)\xrightarrow{\text{a.s.}}K(\theta)$.

PROOF. Since $f_{\theta,m}(y)/f_\theta(y)-1=(\mathbb{P}_m-P)f_\theta(\cdot|y)/h(\cdot)$ by condition (5),

$$\|f_{\theta,m}(\cdot)/f_\theta(\cdot)-1\|_\mathcal{Y}\xrightarrow{\text{a.u.}}_m 0$$

by Lemma 1.9.2 in [32]. This implies

(18) $$\sup_{n\in\mathbb{N}}|K_{m,n}(\theta)-K_n(\theta)|\xrightarrow{\text{a.u.}}_m 0$$

since $K_{m,n}(\theta)-K_n(\theta)=\frac{1}{n}\sum_{j=1}^n\log[f_\theta(Y_j)/f_{\theta,m}(Y_j)]$. Since $K_n(\theta)\xrightarrow{\text{a.s.}}_n K(\theta)$ by the strong law of large numbers, the result follows by the triangle inequality. □

LEMMA A.2. *Under conditions* (3) *and* (4) *of Theorem* 2.2,

(19) $$\liminf_{(m,n)\to(\infty,\infty)}\inf_{\phi\in B}K_{m,n}(\phi)\geq-Q\log P\sup_{\phi\in B}\frac{f_\phi(X,Y)}{h(X)g(Y)}$$

*with probability one for each subset $B$ of $B_\theta\cap C_\theta$.*

PROOF. By condition (3) the term on the right-hand side in (19) is not $-\infty$. Next

$$\inf_{\phi \in B} K_{m,n}(\phi) = -\sup_{\phi \in B} \mathbb{Q}_n \log \mathbb{P}_m \frac{f_\phi(X,Y)}{h(X)g(Y)} \geq -\mathbb{Q}_n \log \mathbb{P}_m \sup_{\phi \in B} \frac{f_\phi(X,Y)}{h(X)g(Y)}.$$

By condition (4), for any $\varepsilon_1 > 0$ and $\varepsilon_2 > 0$, there are measurable $A$ and $M \in \mathbb{N}$ such that $\Pr(A) \geq 1 - \varepsilon_1$ and $\mathbb{P}_m \sup_{\phi \in B} f_\phi(\cdot, y)/h(\cdot)g(y) \leq P \sup_{\phi \in B} f_\phi(\cdot, y)/h(\cdot)g(y) + \varepsilon_2$ for all $m \geq M$ and $y \in \mathcal{Y}$ uniformly on $A$. Hence

$$\inf_{\phi \in B} K_{m,n}(\phi) \geq -\mathbb{Q}_n \log\left\{ P \sup_{\phi \in B} \frac{f_\phi(X,Y)}{h(X)g(Y)} + \varepsilon_2 \right\}$$

for all $m \geq M$ and $y \in \mathcal{Y}$ uniformly on $A$. By the strong law of large numbers on the right-hand side, there are measurable $B$ and $N \in \mathbb{N}$ such that $\Pr(B) \geq 1 - \varepsilon_3$ and

$$-\mathbb{Q}_n \log\left\{ P \sup_{\phi \in B} \frac{f_\phi(X,Y)}{h(X)g(Y)} + \varepsilon_2 \right\} \geq -Q \log\left\{ P \sup_{\phi \in B} \frac{f_\phi(X,Y)}{h(X)g(Y)} + \varepsilon_2 \right\} - \varepsilon_4$$

for all $n \geq N$ uniformly on $B$. Hence

$$\inf_{\phi \in B} K_{m,n}(\phi) \geq -Q \log\left\{ P \sup_{\phi \in B} \frac{f_\phi(X,Y)}{h(X)g(Y)} + \varepsilon_2 \right\} - \varepsilon_4$$

for all $m \geq M$ and $n \geq N$ uniformly on $A \cap B$. We are done since the $\varepsilon$'s were arbitrary. □

LEMMA A.3. *Under conditions* (2) *and* (3) *of Theorem* 2.2, $K$ *is lower semicontinuous.*

PROOF. Let $\theta$ be a point of $\Theta$ and $\{\theta_k\}$ a sequence in $\Theta$ converging to $\theta$. Then

$$\limsup_{k \to \infty} Q \log \frac{f_{\theta_k}(\cdot)}{g(\cdot)} \leq \lim_{n \to \infty} Q \log P \sup_{k \geq n} \frac{f_{\theta_k}(X,Y)}{h(X)g(Y)}$$

$$= Q \log P \limsup_{k \to \infty} \frac{f_{\theta_k}(X,Y)}{h(X)g(Y)},$$

where the equality follows from the monotone convergence theorem by condition (3). Also,

$$\liminf_{k \to \infty} K(\theta_k) \geq -Q \log P \limsup_{k \to \infty} \frac{f_{\theta_k}(X,Y)}{h(X)g(Y)} \geq -Q \log P \frac{f_\theta(X,Y)}{h(X)g(Y)} = K(\theta),$$

where the last inequality follows from condition (2). □

**A.2. Proof of Theorem 2.3.** If we define

$$D_{m,n} = \int_0^1 \nabla^2 K_{m,n}(\theta^* + s(\hat{\theta}_{m,n} - \theta^*))\,ds, \tag{20}$$

then by Taylor series expansion

$$\nabla K_{m,n}(\hat{\theta}_{m,n}) - \nabla K_{m,n}(\theta^*) = D_{m,n}(\hat{\theta}_{m,n} - \theta^*).$$

If we show

$$\left(\frac{V}{n} + \frac{W}{m}\right)^{-1/2} \nabla K_{m,n}(\theta^*) \xrightarrow{\mathcal{L}} \mathcal{N}(0, I), \tag{21}$$

then since $D_{m,n} \xrightarrow{P} J$ by Lemma A.6 below, eventually $D_{m,n}^{-1}$ will exist, and by Slutsky's theorem

$$\left(\frac{V}{n} + \frac{W}{m}\right)^{-1/2} J(\hat{\theta}_{m,n} - \theta^*) = -\left(\frac{V}{n} + \frac{W}{m}\right)^{-1/2} JD_{m,n}^{-1}\nabla K_{m,n}(\theta^*) + o_p(1)$$

$$\xrightarrow{\mathcal{L}} \mathcal{N}(0, I).$$

If we prove (21) under the condition $n/(m+n) \to \alpha$, the subsequence principle gives us (21) without this condition. If $0 < \alpha < 1$, then $(m+n)(V/n + W/m) \to V/\alpha + W/(1-\alpha)$. Since

$$\nabla K_{m,n}(\theta^*) = -\mathbb{Q}_n \nabla \log f_{\theta^*}(Y) - \mathbb{Q}_n \nabla \log \mathbb{P}_m f_{\theta^*}(X|Y)/h(X),$$

we have $\sqrt{m+n}\nabla K_{m,n}(\theta^*) \xrightarrow{\mathcal{L}} \mathcal{N}(0, V/\alpha + W/(1-\alpha))$ by Lemma A.4 below, and in turn we have (21) by Slutsky's theorem. The $\alpha = 0$ and $\alpha = 1$ cases are similar.

LEMMA A.4. *Under conditions* (1) *through* (3), (5) *and* (6) *of Theorem* 2.3,

$$\begin{pmatrix} \sqrt{n}\mathbb{Q}_n \nabla \log f_{\theta^*}(Y) \\ \sqrt{m}\mathbb{Q}_n \nabla \log \mathbb{P}_m f_{\theta^*}(X|Y)/h(X) \end{pmatrix} \xrightarrow{\mathcal{L}} \mathcal{N}\left(0, \begin{pmatrix} V & 0 \\ 0 & W \end{pmatrix}\right). \tag{22}$$

PROOF. By condition (5), $\mathbb{P}_m \xrightarrow{\text{a.u.}} P$ in $l^\infty(\mathcal{F}_2)$ and by condition (6), $\mathbb{G}_m \xrightarrow{\mathcal{L}^*} \mathbb{G}_P$ in $l^\infty(\mathcal{F}_3)$, where $\mathbb{G}_m = \sqrt{m}(\mathbb{P}_m - P)$ and $\mathbb{G}_P$ is a tight Gaussian process in $l^\infty(\mathcal{F}_3)$ with zero mean and covariance function $E(\mathbb{G}_P f \cdot \mathbb{G}_P g) = Pfg - PfPg$. By Slutsky's theorem ([32], Example 1.4.7), $(\mathbb{P}_m, \mathbb{G}_m) \xrightarrow{\mathcal{L}^*} (P, \mathbb{G}_P)$ in $\mathbb{D} = l^\infty(\mathcal{F}_2) \times l^\infty(\mathcal{F}_3)$.

By the almost sure representation theorem ([32], Theorem 1.10.4 and Addendum 1.10.5), if $(\Omega, \mathcal{A}, \Pr)$ is the probability space where $\mathbb{P}_n$ are defined (Pr can be $P^\infty$), there are measurable perfect functions $\phi_m$ on some probability space $(\widetilde{\Omega}, \widetilde{\mathcal{A}}, \widetilde{\Pr})$ such that the following diagram commutes

$$\begin{array}{ccc} \Omega & \xrightarrow{(\mathbb{P}_m, \mathbb{G}_m)} & \mathbb{D} \\ \phi_m \uparrow & \nearrow_{(\widetilde{\mathbb{P}}_m, \widetilde{\mathbb{G}}_m)} & \\ \widetilde{\Omega} & & \end{array}$$

and $\Pr = \widetilde{\Pr} \circ \phi_m^{-1}$ and $(\widetilde{\mathbb{P}}_m, \widetilde{\mathbb{G}}_m) \xrightarrow{\text{a.s.}^*} (\widetilde{\mathbb{P}}_\infty, \widetilde{\mathbb{G}}_\infty)$ in $\mathbb{D}$, where $(\mathbb{P}_\infty, \mathbb{G}_\infty) = (P, \mathbb{G}_P)$ and $(\widetilde{\mathbb{P}}_\infty, \widetilde{\mathbb{G}}_\infty) = (\widetilde{P}, \widetilde{\mathbb{G}}_P)$. Hence for almost all $\tilde{\omega}$, $\sup_{y\in\mathcal{Y}} |(\widetilde{\mathbb{P}}_m - \widetilde{P})(\tilde{\omega}) f_{\theta^*}(\cdot|y)/h(\cdot)| \to 0$ and $\sup_{y\in\mathcal{Y}} |(\widetilde{\mathbb{G}}_m - \widetilde{\mathbb{G}}_P)(\tilde{\omega}) \nabla f_{\theta^*}(\cdot|y)/h(\cdot)| \to 0$. By the uniform continuity of $(s,t) \mapsto t/s$ on $[s_0, \infty) \times \mathbb{R}$ with $s_0 > 0$

$$(23) \quad \sup_{y\in\mathcal{Y}} \left| \frac{\widetilde{\mathbb{G}}_m(\tilde{\omega}) \nabla f_{\theta^*}(\cdot|y)/h(\cdot)}{\widetilde{\mathbb{P}}_m(\tilde{\omega}) f_{\theta^*}(\cdot|y)/h(\cdot)} - \widetilde{\mathbb{G}}_P(\tilde{\omega}) \nabla f_{\theta^*}(\cdot|y)/h(\cdot) \right| \to 0.$$

If we define

$$(24) \quad k(\omega, y) = \mathbb{G}_P(\omega) \nabla f_{\theta^*}(\cdot|y)/h(\cdot),$$

and show $y \mapsto k(\omega, y)$ is bounded and continuous for almost all $\omega$, then the second term on the left-hand side in (23) [which equals $\tilde{k}(\tilde{\omega}, \cdot) = k(\phi_\infty(\tilde{\omega}), \cdot)$] is bounded and continuous for almost all $\tilde{\omega}$. By Lemma A.5 below, $\mathbb{Q}_n(\eta) \widetilde{\mathbb{G}}_P(\tilde{\omega}) \nabla f_{\theta^*}(X|Y)/h(X) \to Q\widetilde{\mathbb{G}}_P(\tilde{\omega}) \nabla f_{\theta^*}(X|Y)/h(X)$ for almost all $\eta$ and $\tilde{\omega}$, and this with (23) leads to

$$(25) \quad \sqrt{m} \mathbb{Q}_n(\eta) \nabla \log \widetilde{\mathbb{P}}_m(\tilde{\omega}) \nabla f_{\theta^*}(X|Y)/h(X) \to Q\widetilde{\mathbb{G}}_P(\tilde{\omega}) \nabla f_{\theta^*}(X|Y)/h(X)$$

for almost all $\eta$ and $\tilde{\omega}$. Even though first $m \to \infty$ and then $n \to \infty$, the limit can be shown to be the same (by a triangle inequality) no matter how $m$ and $n$ go to infinity because of the uniformity in (23).

The function $y \mapsto k(\omega, y)$ is bounded since $\sup_{y\in\mathcal{Y}} |k(\omega, y)| = \|\mathbb{G}_P(\omega)\|_{\mathcal{F}_3} < \infty$ from $\mathbb{G}_P(\omega) \in l^\infty(\mathcal{F}_3)$. Every subscript $i$ refers to the $i$th coordinate in $\mathbb{R}^d$. For almost all $\omega$, the sample path $f \mapsto \mathbb{G}_P(\omega) f$ is $\rho$-continuous on $\mathcal{F}_{3i}$ ([32], Section 1.5), where $\rho(f,g) = \{P(f-g)^2\}^{1/2}$. The function $y \mapsto [\nabla f_{\theta^*}(\cdot|y)/h(\cdot)]_i$ from $\mathcal{Y}$ to $(\mathcal{F}_{3i}, \rho)$ is continuous by condition (2) and the dominated convergence theorem applied to $\rho([\nabla f_{\theta^*}(\cdot|y_n)/h(\cdot)]_i, [\nabla f_{\theta^*}(\cdot|y)/h(\cdot)]_i)^2 \leq 4P(F_i^2) < \infty$ with $y_n \to y$ and $F$ in condition (6). The function $y \mapsto k_i(\omega, y)$ is a composition of the two continuous functions, hence continuous, for almost all $\omega$.

By the central limit theorem $\sqrt{n} \mathbb{Q}_n \nabla \log f_{\theta^*}(Y) \xrightarrow{\mathcal{L}} \mathcal{N}(0, V)$, and if $(H, \mathcal{B}, \text{Qr})$ is the probability space where $\mathbb{Q}_n$ are defined (Qr can be $Q^\infty$), there is an almost sure representation for this with commutative diagram

$$\begin{array}{ccc} H & \xrightarrow{\mathbb{Q}_n} & \mathbb{R}^d \\ \psi_m \uparrow & \nearrow & \\ \widetilde{H} & \widetilde{\mathbb{Q}}_n & \end{array}$$

and $\mathrm{Qr} = \widetilde{\mathrm{Qr}} \circ \psi_n^{-1}$. If we combine this representation with (25),

$$\begin{pmatrix} \sqrt{n}\widetilde{\mathbb{Q}}_n(\tilde{\eta})\nabla \log f_{\theta^*}(Y) \\ \sqrt{m}\widetilde{\mathbb{Q}}_n(\tilde{\eta})\nabla \log \widetilde{\mathbb{P}}_m(\tilde{\omega})\nabla f_{\theta^*}(X|Y)/h(X) \end{pmatrix}$$
$$\to_{m,n} \begin{pmatrix} Z(\tilde{\eta}) \\ Q\widetilde{\mathbb{G}}_P(\tilde{\omega})\nabla f_{\theta^*}(X|Y)/h(X) \end{pmatrix}$$

for almost all $\tilde{\eta}$ and $\tilde{\omega}$, where $Z(\tilde{\eta})$ is $\mathcal{N}(0,V)$. In this representation, it is clear that the two terms on the right-hand side, being functions of independent random variables, are independent. This almost sure convergence implies weak convergence, and undoing the almost sure representation gives

$$\begin{pmatrix} \sqrt{n}\mathbb{Q}_n\nabla \log f_{\theta^*}(Y) \\ \sqrt{m}\mathbb{Q}_n\nabla \log \mathbb{P}_m\nabla f_{\theta^*}(X|Y)/h(X) \end{pmatrix} \xrightarrow{\mathcal{L}}_{m,n} \begin{pmatrix} Z \\ Q\mathbb{G}_P\nabla f_{\theta^*}(X|Y)/h(X) \end{pmatrix}.$$

We are done if we show that the second term on the right-hand side is $\mathcal{N}(0,W)$.

Let $T = Q\,\mathbb{G}_P \nabla f_{\theta^*}(X|Y)/h(X)$. Note $T(\omega) = Qk(\omega,\cdot)$ with $k$ in (24). By condition (2) there is a sequence $\{Q_i\}$ of probability measures with finite support such that $Q_i \xrightarrow{\mathcal{L}} Q$ ([1], Theorem 14.10 and Theorem 14.12). Let $T_i(\omega) = Q_i k(\omega,\cdot)$. Then $T_i(\omega) \to T(\omega)$ for almost all $\omega$ because $y \mapsto k(\omega,y)$ is bounded and continuous for almost all $\omega$. Since $\mathbb{G}_P$ is a Gaussian process, $T_i$ is normally distributed. By condition (6), $(y,s) \mapsto E[k(\cdot,y)k(\cdot,s)^T]$ is bounded and continuous by the dominated convergence theorem. Hence $\mathrm{var}\,T_i \to \mathrm{var}\,T$, and by Fubini the limit equals $W$. Now for any $t \in \mathbb{R}^d$ $\exp(-t^T(\mathrm{var}\,T_i)t/2) \to \exp(-t^T W t/2)$. Hence $T_i \xrightarrow{\mathcal{L}} \mathcal{N}(0,W)$ and $T \sim \mathcal{N}(0,W)$.
□

LEMMA A.5. *Under condition* (2) *of Theorem* 2.3, $\mathbb{Q}_n \xrightarrow{\mathcal{L}} Q$ *almost surely.*

PROOF. Let $\mathcal{B}$ be a countable basis for $\mathcal{Y}$ and $\mathcal{A}$ be the set of all finite intersections of elements of $\mathcal{B}$ (also countable). For each $A \in \mathcal{A}$ we have $\mathbb{Q}_n(A) \to Q(A)$ by the strong law of large numbers. Hence, a countable union of null sets being a null set, this holds simultaneously for all $A \in \mathcal{A}$. The result follows since $\mathcal{A}$ is a convergence determining class ([4], Theorem 2.2).
□

LEMMA A.6. *Under conditions* (4) *through* (8) *of Theorem* 2.3, $D_{m,n} \xrightarrow{P} J$, *where* $D_{m,n}$ *is defined by* (20) *and* $J$ *in condition* (3) *of Theorem* 2.3.

PROOF. First note

$$|D_{m,n} - J| \leq \int_0^1 |\nabla^2 K_{m,n}(\theta^* + s(\hat{\theta}_{m,n} - \theta^*)) + Q\nabla^2 \log f_{\theta^* + s(\hat{\theta}_{m,n} - \theta^*)}(Y)| \, ds$$

$$+ \sup_{0 \leq s \leq 1} |Q\nabla^2 \log f_{\theta^* + s(\hat{\theta}_{m,n} - \theta^*)}(Y) - Q\nabla^2 \log f_{\theta^*}(Y)|.$$

By condition (4), $Q\nabla^2 \log f_\theta(Y)$ is continuous on $S_\rho$ ([9], page 110). Hence the second term on the right-hand side converges in probability to zero by the weak consistency of $\hat{\theta}_{m,n}$. The first term on the right-hand side will also converge in probability to zero because for any $\varepsilon > 0$

$$(26) \quad \Pr\left(\int_0^1 |\nabla^2 K_{m,n}(\theta^* + s(\hat{\theta}_{m,n} - \theta^*)) + Q\nabla^2 \log f_{\theta^* + s(\hat{\theta}_{m,n} - \theta^*)}(Y)| \, ds > \varepsilon\right)$$

$$\leq \Pr(\hat{\theta}_{m,n} \notin S_\rho) + \Pr\left(\sup_{\theta \in S_\rho} |\nabla^2 K_{m,n}(\theta) + Q\nabla^2 \log f_\theta(Y)| > \varepsilon\right),$$

if we show the second term on the right-hand goes to zero. Note $\nabla^2 K_{m,n}(\theta) = -\mathbb{Q}_n \nabla^2 \log f_\theta(Y) - \mathbb{Q}_n W_m(\theta, Y)$ where

$$W_m(\theta, y) = \frac{\mathbb{P}_m \nabla^2 f_\theta(\cdot|y)/h(\cdot)}{\mathbb{P}_m f_\theta(\cdot|y)/h(\cdot)} - \frac{\{\mathbb{P}_m \nabla f_\theta(\cdot|y)/h(\cdot)\}\{\mathbb{P}_m \nabla f_\theta(\cdot|y)/h(\cdot)\}^T}{\{\mathbb{P}_m f_\theta(\cdot|y)/h(\cdot)\}^2}.$$

By condition (4), $\sup_{\theta \in S_\rho} |\mathbb{Q}_n \nabla^2 \log f_\theta(Y) - Q\nabla^2 \log f_\theta(Y)| \xrightarrow{\text{a.s.}^*} 0$. Hence the second term on the right-hand side in (26) will go to zero, if we show $\sup_{\theta \in S_\rho} |\mathbb{Q}_n W_m(\theta, Y)| \xrightarrow{\text{a.s.}^*} 0$.

By condition (7)

$$(27) \quad \sup_{\theta \in S_\rho} \sup_{y \in \mathcal{Y}} |\mathbb{P}_m \nabla^2 f_\theta(\cdot|y)/h(\cdot)| \xrightarrow{\text{a.s.}^*} 0.$$

Expanding $\mathbb{P}_m \nabla f_\theta(\cdot|y)/h(\cdot)$ as

$$\mathbb{P}_m \nabla f_\theta(\cdot|y)/h(\cdot) = \mathbb{P}_m \nabla f_{\theta^*}(\cdot|y)/h(\cdot)$$

$$+ \int_0^1 \mathbb{P}_m \nabla^2 f_{\theta^* + s(\theta - \theta^*)}(\cdot|y)/h(\cdot)(\theta - \theta^*) \, ds$$

leads to, for any $\theta \in S_\rho$,

$$\sup_{y \in \mathcal{Y}} |\mathbb{P}_m \nabla f_\theta(\cdot|y)/h(\cdot)| \leq \sup_{y \in \mathcal{Y}} |\mathbb{P}_m \nabla f_{\theta^*}(\cdot|y)/h(\cdot)|$$

$$+ \sup_{\theta \in S_\rho} \sup_{y \in \mathcal{Y}} |\mathbb{P}_m \nabla^2 f_\theta(\cdot|y)/h(\cdot)| \rho.$$

The first term on the right-hand side converges almost surely to zero because $\mathcal{F}_3$ is $P$-Glivenko–Cantelli from being $P$-Donsker [condition (6)]. Since the second term on the right-hand side also converges almost surely to zero by (27),

$$\sup_{\theta \in S_\rho} \sup_{y \in \mathcal{Y}} |\mathbb{P}_m \nabla f_\theta(\cdot|y)/h(\cdot)| \xrightarrow{\text{a.s.}^*} 0. \tag{28}$$

With condition (5) and (28), $\sup_{\theta \in S_\rho} \sup_{y \in \mathcal{Y}} |\mathbb{P}_m f_\theta(\cdot|y)/h(\cdot) - 1| \xrightarrow{\text{a.s.}^*} 0$, and this, (27) and (28) imply $\sup_{\theta \in S_\rho} \sup_{y \in \mathcal{Y}} |W_m(\theta, y)| \xrightarrow{\text{a.s.}^*} 0$. We are done because $\sup_{\theta \in S_\rho} |\mathbb{Q}_n W_m(\theta, Y)| \leq \sup_{\theta \in S_\rho} \sup_{y \in \mathcal{Y}} |W_m(\theta, y)|$. □

**Acknowledgments.** We thank two anonymous reviewers, an Associate Editor and the Editor for advice on improving the manuscript.

DIVISION OF BIOSTATISTICS
SCHOOL OF MEDICINE
WASHINGTON UNIVERSITY
ST. LOUIS, MISSOURI 63110-1093
USA
E-MAIL: yunju@wubios.wustl.edu

SCHOOL OF STATISTICS
UNIVERSITY OF MINNESOTA
313 FORD HALL
224 CHURCH ST. S. E.
MINNEAPOLIS, MINNESOTA 55455
USA
E-MAIL: charlie@stat.umn.edu